\documentclass[11pt]{amsart}
\usepackage{enumerate}
\newtheorem{lema}{Lemma}[section]
\newtheorem{theo}[lema]{Theorem}

\newtheorem{coro}[lema]{Corollary}

\newtheorem{rema}[lema]{Remark}
\usepackage{graphicx}
\title[]
{Polynomial differential equations with small coefficients}
\author[]
{M. A. M. Alwash}
\address{Departments of Mathematics, West Los Angeles College and University of California, Los Angeles
\newline
9000 Overland Ave, Culver City, CA 90230-3519, USA}
\email{alwashm@wlac.edu
\newline
April 2, 2008}
\begin{document}
\maketitle
\numberwithin{equation}{section}
\begin{abstract}
 Classes of polynomial non-autonomous differential equations of degree $n$ are considered. An explicit bound on the size
 of the coefficients is given which implies that each equation in the class has exactly $n$ complex
 periodic solutions. In most of the classes the upper bound can be improved when we consider
  real periodic solutions.
 We present a proof to a recent conjecture about the number of periodic solutions.
 The results are used to give upper bounds for the number of limit cycles of polynomial two-dimensional systems.
 \vspace{0.25 cm}
 \vspace{0.25cm}
\end{abstract}
\section{Introduction}
We consider differential equations of the form
\begin{equation}
\label{e1.1}
\dot{z}:=\frac{dz}{dt}=z^{n}+P_{n-1}(t)z^{n-1}+\dots+ P_{1}(t)z+P_{0}(t)
\end{equation}
where $z$ is a complex-valued function and $P_{i}$ are real-valued continuous functions. This class of equations has received some
attention in the literature. The main concern is to estimate the number of periodic solutions. The qualitative behavior of the
solution curves depends entirely on the periodic solutions; see, for example, \cite{p3}. The problem was suggested by C. Pugh
as a version of Hilbert's sixteenth problem; it is listed as Problem 7 by Steve Smale in \cite{s1}.

We denote by $z(t,c)$ the solution of \eqref{e1.1} satisfying $z(0,c)=c$. For a fixed real number $\omega$, we define
the set $Q$ to be the set of all complex numbers $c$ such that $z(t,c)$ is defined for all $t$ in the interval
$[0,\omega]$; the set $Q$ is an open set. On $Q$ we define the displacement function $q$ by
\[
q(c)=z(\omega,c)-c.
\]
Zeros of $q$ identify initial points of solutions of \eqref{e1.1} which satisfy the boundary condition
$z(0)=z(\omega)$. We describe such solutions as {\em periodic} even when the functions $P_{i}$ are not themselves
periodic. However, if $P_{i}$ are $\omega$-periodic then these solutions are also $\omega$-periodic.
\newpage
Note that $q$ is holomorphic on $Q$. The multiplicity of a periodic solution $\varphi$ is that of $\varphi (0)$ as a
zero of $q$. It is useful to work with a complex dependent variable. The reason is that periodic solutions cannot then
be destroyed by small perturbations of the right-hand side of the equation. Suppose that $\varphi$ is a periodic
solution of multiplicity $k$. This solution is counted as $k$ solutions.
By applying Rouche's theorem to the function $q$, for any sufficiently small
perturbations of the equation, there are precisely $k$ periodic solutions in a neighborhood of $\varphi$ (counting
multiplicity).

Upper bounds to the number of periodic solutions of \eqref{e1.1} can be used as upper bounds to the
number of periodic solutions when $z$ is limited to be real-valued.
This is the reason that $P_{i}$ are not allowed to
be complex-valued. It was shown in \cite{a5} that when the coefficients are complex, then there are centers. We recall that
a \emph{center} is a family of periodic solutions with an open set of initial conditions.

When $n=3$, equation \eqref{e1.1} is known as the Abel differential equation. This case is of particular interest
because of a connection with Hilbert's sixteenth problem. It was shown in \cite{l2} and
\cite{p3} that the Abel differential equation has exactly three periodic solutions provided account
is taken of multiplicity. Local questions related to Hilbert's sixteenth problem (bifurcation of
small-amplitude limit cycles and center conditions) are reduced to polynomial equations
in which the leading coefficient
changes sign. However, we show in Section 3 that the result of this paper can provide upper bounds for the number of limit cycles.

Lins Neto \cite{l1} has given examples
which demonstrate that there is no upper bound for the number of periodic solutions for cubic equations when the
leading coefficient changes sign. Panov \cite{p2} demonstrated that there is no upper for the number of periodic
solutions by proving that the set of return maps of these equations are dense in the space of orientation preserving homeomorphisms.
It follows from these examples that there is no upper bound, in terms of $n$ only, for the number of periodic solutions of equation (1.1).

The case $n=4$ was considered in \cite{a3} and \cite{a4}. The main concern was the multiplicity of
$z=0$ when the coefficients are polynomial functions in $t$, and in $\cos{t}$ and $\sin{t}$. Equations with
at least $10$ real periodic solutions were constructed. These periodic solutions are bifurcated from a periodic solution of
multiplicity $10$.

Our aim in this paper is to gain information on the total number of periodic solutions;
this is a global question, while looking at multiplicity leads only to local results.
Ilyashenko \cite{i1}, proved that the number of real periodic solutions does not exceed
$ 8\exp[(3C+2)\exp(1.5(2C+3)^{n})], $
where $C>1$ is an upper bound for $|P_{i}(t)|$.
Although this bound is non-realistic, it is the only known explicit estimate. An open problem was made in \cite{i1} to find a similar upper bound for
complex differential equations; that is for complex solutions and complex coefficients. However, such equations could have centers. We
conjectured in \cite{a5} that equations with polynomial coefficients do not have centers when $P_{0}\equiv 0$ and $P_{1}\equiv 0$

Classes of equations that contain only three terms were considered in \cite{a1}, \cite{a2}, \cite{g1}, and \cite{p1}.
Under certain conditions, equations in these classes have exactly $n$ complex periodic solutions.
Now, we consider the general case and without a restriction on the number of terms. A {\em singular periodic solution} is the limit of a sequence
of periodic solutions which is undefined on the interval $[0,\omega]$. If the equations in a class do not have singular periodic solutions then all the
equations in a path component have the same number of periodic solutions. We present two classes such that
equations in each class do not have singular periodic solutions. It is shown that the qualitative behavior of the solutions is dominated by the term
$z^{n}$ when $|z|>1$, and by the linear term, $P_{1}(t)z$ when $|z|<1$.
In particular, we prove the following results:

\begin{theo}
Let $K$ be a real number such that $K>\max{|P_{0}(t)|}$.
If for $0 \leq t \leq \omega$
 \[ |P_{i}(t)|  \leq \frac{K-|P_{0}(t)|}{(n-3)\sqrt[n]{K^{i}}}\sin{(\frac{\pi}{n-2})},i \neq 0, 1, n-1\]
 and
\[|P_{1}(t)| \geq \frac{K+|P_{0}(t)|+|P_{n-1}(t)|+(K-|P_{0}(t)|)\sin{(\frac{\pi}{n-2})}}{\sqrt[n]{K}}\]
 then  \eqref{e1.1} has $n$ complex periodic solutions. If $P_{n-1}(t)$ satisfies the above condition, then there are
 at most two positive periodic solutions and
 at most two negative periodic solutions.
 In the case $n$ is odd, there is at least one real periodic solution and there are at most three real periodic solutions.
\end{theo}

\begin{theo}
Let $K$ be a real number such that
$K>\max{|P_{0}(t)|}$.
If for $0 \leq t \leq \omega$
 \[P_{n-2}(t)\leq 0,\]
 \[ |P_{i}(t)| \leq \frac{(K-\sqrt[n]{K^{n-2}}P_{n-2}(t)-|P_{0}(t)|)}{(n-4)\sqrt[n]{K^{i}}}\sin{(\frac{\pi}{n-2})},\,\,
 i \neq 0, 1, n-2\]
 and
\[|P_{1}(t)| \geq \frac{K+|P_{0}(t)|-\sqrt[n]{K^{n-2}}P_{n-2}(t)+\frac{n-3}{n-4}(K-\sqrt[n]{K^{n-2}}P_{n-2}(t)-|P_{0}(t)|)
\sin{(\frac{\pi}{n-2})}}{\sqrt[n]{K}}\]
 then  \eqref{e1.1} has exactly $n$ complex periodic solutions.
\end{theo}

Since the coefficients $P_{i}(t)$ are real functions, it follows from Theorems 1.1 and 1.2 that equation (1.1) has at
most $n$ real periodic solutions. These results generalize a similar result of Calanchi and Ruf \cite{c1}.
It was shown in \cite{c1} that if $n$ is odd and
$|P_{i}(t)| \leq \frac{1}{2n(n-1)}, i=1,2,\cdots,n-1$, then
the equation has at most $n$ real periodic solutions. It was conjectured in \cite{c2} that the oddness assumption on the
degree, $n$, is not necessary. Theorem 1.1 shows that the conjecture is true if $|P_{1}(t)|$ is large.
Our proof is much shorter than that of \cite{c1} and $n$ is not restricted to be odd. Moreover, our upper bound for the number of
real periodic solutions is constant and does
not depend on $n$. In our results and the result of \cite{c1}, the size of small coefficients is of order $\frac{1}{n^{2}}$.

Panov \cite{p1} proved that the equation has $n$ periodic solutions
of period $\frac{1}{6n}$ when the roots of the right hand side in \eqref{e1.1} are inside the disc $|z| < \frac{\pi}{2n}$.
We explain in Section 3 that the result of \cite{p1} is different from the results given in this paper. The condition of
\cite{p1} implies that all periodic solutions enter the unit disk; in our case, periodic solutions lie outside the unit disk.
It was shown in \cite{l2} that there exists an $\varepsilon > 0$ such that if $|P_{i}(t)| < \varepsilon$ for $i \neq 1$,
then the equation \eqref{e1.1} has $n$ complex periodic solutions. The same conclusion holds when
$n$ is odd and $|P_{i}(t)| < \varepsilon$ for $i \neq 1, \frac{1}{2}(n+1)$.

Now, we give upper bounds for the number of real periodic solutions when $|P_{2}(t)|$, $|P_{1}(t)|$, or $|P_{0}(t)|$ is large.
\begin{theo}
\begin{itemize}
\item[(i)]
If for $0 \leq t \leq \omega$
\[|P_{i}(t)| < \frac{n}{(n-1)^{2}},\,\,1 \leq i \leq n-1,\]
and
\[|P_{0}(t)| > \frac{2n-1}{n-1} \]
then \eqref{e1.1} has at most one positive real periodic solution and at most one real negative periodic solution.
\item[(ii)]
If for $0 \leq t \leq \omega$
\[|P_{i}(t)| < \frac{n}{(n-2)^{2}},\,\,2 \leq i \leq n-1,\]
and
\[|P_{1}(t)| > \frac{n(2n-3)}{n-2} \]
then \eqref{e1.1} has at most five real periodic solutions; at most three of them are positive and at most three
of them are negative.
\item[(iii)]
If for $0 \leq t \leq \omega$
\[|P_{i}(t)| < \frac{n}{(n-3)^{2}},\,\,3 \leq i \leq n-1,\]
and
\[|P_{2}(t)| > \frac{n(n-1)(2n-5)}{2(n-3)} \]
then \eqref{e1.1} has at most eight real periodic solutions; at most five of them are positive and at most five of
them are negative.
\end{itemize}
\end{theo}
In the second part of Theorem 1.3, no condition is imposed on $P_{0}$, and in the third part no condition is imposed on
$P_{0}$ and $P_{1}$.

If the leading coefficient, $P_{n}(t)$, in equation (1.1) is not $1$ but does not vanish anywhere then the transformation of
the independent variable reduces the equation into a similar equation but with a leading coefficient equals 1.
Bounds similar to those in the above results are obtained in Section 3.

In Section 2, we describe the
phase portrait of \eqref{e1.1} and recall some results from \cite{l2} and \cite{l3}. In Section 3, we present the proofs of our results.
We also demonstrate in Section 4 that our results can be used to give upper bounds for the number of limit cycles for certain classes of
polynomial two-dimensional systems.

\section{The Phase Portrait}
We identify equation \eqref{e1.1} with the vector $(P_{0},P_{1},P_{2},\cdots,P_{n-1})$ and write $\mathcal{L}$ for the set of
all equations of this form. With the usual definitions of additions and scalar multiplications, $\mathcal{L}$ is a
linear space; it is a normed space if for
 $P=(P_{0},P_{1},P_{2},\cdots,P_{n-1})$, we define
 \[
 \| P \| = \max\{\max_{0\leq t \leq \omega}|P_{0}(t)|,
\max_{0\leq t \leq \omega}|P_{1}(t)|, \max_{0\leq t \leq \omega}|P_{2}(t)|, \cdots,\max_{0\leq t \leq \omega}|P_{n-1}(t)| \}
\]
The displacement function $q$ is holomorphic on the open set $Q$.
Moreover, $q$ depends continuously on $P$ with the above norm on $\mathcal{L}$ and the topology of uniform convergence
on compact sets on the set of holomorphic functions.
 If $\varphi$ is a non-real solution which is
periodic, then so is $\bar{\varphi}$, its complex conjugate.

In \cite{l2}, it was shown that the phase portrait of \eqref{e1.1} is as shown in Figure 1 below. We refer to \cite{l2}
for the details. There, the coefficients $P_{i}(t)$ were $\omega-$periodic. It can be verified that the same methods are
applicable to the study of the number solutions that satisfy $z(0)=z(\omega)$ whether the coefficients are periodic or
not.

Note that the radius, $\rho$, of the disc $D$ depends only on $\|P\|$ and $\omega$. For any given equation, $\rho$ can be
determined explicitly. The disc $D$ is the larger one in the figure.
 If $z=r e^{i\theta}$ then the sets
$G_{k},k=0,1,\dots,2n-3$, which are the arms in the figure, are defined by
\[
G_{k}= \{z| r > \rho, \frac{k\pi}{n-1}-\frac{a}{r} < \theta < \frac{k\pi}{n-1}+\frac{a}{r}\}
\]
where $a=\max\{6,6\|P\|\}$. Between the arms are the sets $H_{k},k=0,1,\dots,2n-3$, which are defined by
\[
H_{k}= \{z| > \rho, \frac{k\pi}{n-1}+\frac{a}{r} \leq \theta \leq \frac{(k+1)\pi}{n-1} -\frac{a}{r}\}
\]

In $G_{k}$, $\dot{r}>0$ when $k$ is even, and $\dot{r}<0$ when $k$ is odd.
In $H_{k}$, $\dot{\theta}>0$ when $k$ is even, and $\dot{\theta}<0$ when $k$ is odd.
For even $k$, trajectories can enter $G_{k}$ only across $r=\rho$, and for odd $k$, trajectories can leave $G_{k}$ only
across $r=\rho$. No solution can become infinite in $H_{k}$ as time either increases or decreases. Every solution enters
$D$. Solutions become unbounded if and only if they remain in one of the arms $G_{k}$, tending to infinity
 as $t$ increases if $k$ is even and as $t$ decreases if $k$ is odd.

\begin{figure}[ht]
\begin{picture}(160,160)
\put(75,75){\circle{20}}
\put(75,75){\circle{55}}
\put(116,76){\vector(0,1){10}}
\put(112,74){\vector(0,-1){10}}
\put(95,78){\line(1,0){45}}
\put(95,72){\line(1,0){45}}
\put(98,74){\vector(-1,0){10}}
\put(130,74){\vector(-1,0){10}}
\put(84,86){\vector(1,-1){10}}
\put(104,96){\vector(1,-1){10}}
\put(129,105){\vector(1,-1){10}}
\put(129,35){\vector(1,1){10}}
\put(98,55){\vector(1,1){10}}
\put(83,61){\vector(1,1){10}}
 \put(85,70){\line(2,-1){37}}
\put(85,78){\line(2,1){37}}
\put(131,93){\makebox(0,0)[b1]{$\theta_{1}$}}
\put(130,42){\makebox(0,0)[b1]{$\theta_{2}$}}
\put(84,94){\line(1,1){35}}\put(89,90){\line(1,1){35}}
\put(85,90){\vector(1,1){10}}\put(105,110){\vector(1,1){10}}
\put(85,59){\line(1,-1){35}}\put(79,56){\line(1,-1){35}}
\put(82,58){\vector(1,-1){10}}\put(103,37){\vector(1,-1){10}}
\put(145,95){\makebox(0,,0)[bl]{$H_{k}$}}
\put(145,70){\makebox(0,,0)[bl]{$G_{k}$}}
\put(145,40){\makebox(0,,0)[bl]{$H_{k-1}$}}
\put(127,124){\makebox(0,0)[bl]{$G_{k+1}$}}
\put(122,10){\makebox(0,0)[bl]{$G_{k-1}$}}
\thicklines
\end{picture}
\caption{Phase Portrait Around $G_{k}$}
\end{figure}

Let $q(P,c)=z_{P}(\omega,c)-c$, where $z_{P}(t,c)$ is the solution of
 $P \in \mathcal{L}$ satisfying $z_{P}(0,c)=c$. Suppose that $(P_{j})$ and
$(c_{j})$ are sequences in $\mathcal{L}$ and $\mathbb{C}$, respectively, such that $q(P_{j},c_{j})=0$. If $P_{j} \to  P$
and $c_{j} \to  c$ as $j \to  \infty$, then either $q(P,c)=0$, in this case $z_{P}(t,c)$ is a periodic solution, or
$z_{P}(t,c)$ is not defined for the whole interval $0 \leq t \leq \omega$. In the later case, we say that $z_{P}(t,c)$
is a {\em singular periodic solution}. We also say that $P$ has a singular periodic solution if $c_{j} \to  \infty$; in
this case there are $\tau$ and $c$ such that the solution $z_{P}$ with $z_{P}(\tau)=c$ becomes unbounded at finite time
as $t$ increases and as $t$ decreases. We summarize the results of \cite{l2} related to the work in this paper.

\begin{lema}
\begin{itemize}
\item[(i)] Let $\mathcal{A}$ be the subset of $\mathcal{L}$ consisting of all equations which have no singular
    periodic solutions. The set $\mathcal{A}$ is open in $\mathcal{L}$. All equations in the same components of
    $\mathcal{A}$ have the same number of periodic solutions.
\item[(ii)] The equation $\dot{z}=z^{n}$ has exactly $n$ periodic solutions.
\end{itemize}
\end{lema}
Now, we give the derivatives of the displacement function $q(c)$, see \cite{l3}.
\begin{lema}
Consider the real differential equation
\begin{equation}
\label{e2.1}
\dot{x}=f(x,t)=x^{n}+P_{n-1}(t)x^{n-1}+\dots+ P_{1}(t)x+P_{0}(t).
\end{equation}
The derivatives of the displacement function are given by
\[ q'(c) = E(\omega,c)-1, \]
\[ q''(c)=E(\omega,c) \int_{0}^{\omega}D(t,c) dt,\]
\[ q'''(c)= E(\omega,c)[\frac{3}{2}(G(\omega,c))^{2}+\int_{0}^{\omega}(E(\omega,c))^{2}f_{3}(x(t,c),t)dt] \]
where
\[ E(t,c) = \exp{ \int_{0}^{t} f_{1}(x(\tau,c),\tau) d\tau}, \]
\[ D(t,c)= E(t,c) f_{2}(x(t,c),t),\]
and
\[ G(t,c)=\int_{0}^{t} D(\tau,c) d\tau.\]
Here, $f_{i}=\frac{\partial^{i}}{\partial x^{i}}$.
\end{lema}
These formulae imply that for $k \leq 3$, if $f_{k}$ does not change sign on an interval then the $k$th derivative of $q$ does not change sign
on that interval.

\section{Proofs of Main Results}
First, we give the proof of Theorem 1.1.
\begin{proof}
\emph{ (Theorem 1.1:)}\\
We make the transformations
\[z \mapsto \frac{1}{\sqrt[n]{K}}\,\,z ,\]
\[ t \mapsto \sqrt[n]{K^{n-1}}\,\, t,\]
the equation becomes
\[ \dot{z}= z^{n}+\frac{P_{n-1}(t)}{\sqrt[n]{K}}z^{n-1}+\frac{P_{n-2}(t)}{\sqrt[n]{K^{2}}}z^{n-2}+\dots+
\frac{P_{1}(t)}{\sqrt[n]{K^{n-1}}}z+
\frac{P_{0}(t)}{K}.\]
These transformations preserve the number of periodic solutions. In this equation the condition
\[|\frac{P_{0}(t)}{K}| < 1\]
is satisfied. Therefore, without lost in the generality, we assume that $K=1$.

 With $z=re^{i\theta}$, we have
\[
r\dot{\theta}=r^{n}\sin{((n-1)\theta)}+r^{n-1}P_{n-1} \sin{((n-2)\theta)}+ \cdots+r^{2}P_{2} \sin{(\theta)}-P_{0}\sin{(\theta)}.
\]
\[
\dot{r}=r^{n}\cos{((n-1)\theta)}+r^{n-1}P_{n-1} \cos{((n-2)\theta)}+ \cdots+r^{2}P_{2} \cos{(\theta)}+P_{1}r+P_{0}\cos{(\theta)}.
\]
 A singular periodic solution enters $D$ from $G_{k}$ with odd $k$ and leaves $D$ to a $G_{j}$ with
 even $j$. We refer the reader to Figure 1, where the large circle is $r=\rho$ and the small circle is $r=1$.

Now, if $k$ is odd and $1 \leq k \leq n-1$, let $\theta_{1}=\frac{k\pi}{n-2}, \theta_{2}=\frac{(k-1)\pi}{n-2}$.
We take $\rho$ sufficiently large such that the half lines $\theta = \theta_{1}$ and $\theta=\theta_{2}$ do not pass
through $G_{k-1},G_{k}$, or $G_{k+1}$. Precisely, we take
\[ \rho > \frac{a(n-1)(n-2)}{\pi}.\]

Consider
\[
r\dot{\theta}(\theta_{1})=r^{n}\sin{((n-1)\theta_{1})}+r^{n-1}P_{n-1} \sin{((n-2)\theta_{1})}+
\cdots+r^{2}P_{2} \sin{(\theta_{1})}- P_{0}\sin{(\theta_{1})}.
\]
Now, we use the properties
\[
\sin{((n-1)\theta_{1})}=-\sin{(\frac{k\pi}{n-2})},\]
\[
\sin{((n-1)\theta_{2})}=\sin{(\frac{(k-1)\pi}{n-2})},\]
\[\sin{((n-2)\theta_{1})}=0,\]
\[\sin{((n-3)\theta_{1})}=\sin{(\frac{k\pi}{n-2})},\]
and
\[\sin{(\frac{k\pi}{n-2})} \geq \sin{(\frac{\pi}{n-2})}.\]

With these facts and for $1<r<\rho$, we have
\[
r\dot{\theta}(\theta_{1}) \leq - r^{n}\sin{(\frac{k\pi}{n-2})}+(n-2)C r^{n-2}+|P_{0}|\sin{(\frac{k\pi}{n-2})}.
\]
Hence, $r\dot{\theta}(\theta_{1}) < 0$ when $r>1$ and $ (n-2)C \leq (1-|P_{0}|)\sin{(\frac{\pi}{n-2})}$,
where, $C$ is the upper bound for $|P_{i}(t)|,\,\,i \neq 0, 1$.

Similarly, we have
\[
r\dot{\theta}(\theta_{2}) \geq r^{n}\sin{(\frac{(k-1)\pi}{n-2})}-(n-2)Cr^{n-2} -|P_{0}|\sin{(\frac{(k-1) \pi}{n-2})}> 0.
\]
Therefore, no singular periodic solution can enter $D$ form $G_{k}$ and leave $D$ when $r>1$. The other possibility is that
such a solution enters and then leaves the disk $r <1$. However,
\[
\dot{r}(r=1)=\cos{((n-1)\theta)}+P_{n-1} \cos{((n-2)\theta)}+ \cdots+P_{2} \cos{(\theta)}+P_{1}+P_{0}\cos{(\theta)}.
\]
With the above conditions,
\[
\dot{r}(r=1) \leq (n-3)C +1+|P_{n-1}(t)|+P_{1} +|P_{0}| \leq 0,
\]
when $P_{1} \leq -(1+|P_{0}| + |P_{n-1}(t)|)- (1-|P_{0}|) \sin{(\frac{\pi}{n-2})}$. In this case, no solution leaves the disk $r<1$.

On the other hand,
\[
\dot{r}(r=1) \geq - (n-3)C - 1 + P_{1} -|P_{0}|-|P_{n-1}(t)| \geq 0,
\]
when $P_{1} \geq 1 + |P_{0}|+|P_{n-1}(t)|+(1-|P_{0}|) \sin{(\frac{\pi}{n-2})}$. In this case, no solution enters the disk $r<1$.

Under the hypotheses in the statement of Theorem 1.1, no singular periodic solution can enter $D$ from $G_{k}$ and leave $D$,
when $k$ is odd and $k \leq n-1$. Since the phase portrait is symmetric about the $x-$axis, if $k$ is odd and $n \leq k \leq 2n-3$,
no singular periodic solution enters $D$ from $G_{k}$ and leaves $D$. It follows that the equation does not
have a singular periodic solution.

Now, consider the class of equations
\[
\dot{z}=z^{n}+\lambda[P_{n-1}(t)z^{n-1}+\cdots+P_{2}(t)z^{2}+P_{0}(t)]+P_{1}(t)z,
\]
with $0 \leq \lambda \leq 1$. For any equation in this family, the conditions in the statement of Theorem 1.1 become
\[\lambda (n-3)C \leq (1-\lambda |P_{0}|)\sin{(\frac{\pi}{n-2})},\]
\[|P_{1}| \geq 1+ \lambda |P_{0}|  +\lambda |P_{n-1}|+(1-\lambda |P_{0}|)\sin{(\frac{\pi}{n-2})}.\]
Dividing the first inequality by $\lambda$ and rearrange the terms in the second inequality give
\[ (n-3)C \leq (\frac{1}{\lambda}- |P_{0}|)\sin{(\frac{\pi}{n-2})},\]
\[|P_{1}| \geq 1+ \lambda |P_{0}|(1-\sin{(\frac{\pi}{n-2}))}+\sin{(\frac{\pi}{n-2})}+\lambda |P_{n-1}|.\]
It is clear that both inequalities are true for any $\lambda \leq 1$ if they are true for $\lambda = 1$. Therefore, any
equation in this class satisfies the hypotheses of Theorem 1.1, and hence does not have
a singular periodic solution. Since the equations in this family are in the same component of $\mathcal{A}$,
it follows from Lemma 2.1, that these equations have the same number of periodic solutions.
The equation
\[ \dot{z}=z^{n}+P_{1}(t)z \]
belongs to this family. This equation has $n$ periodic solutions. To prove this statement, we follow the same procedure. In this case
\[
r \dot{\theta}(\theta_{1})=r^{n} \sin{((n-1)\theta_{1})} \leq 0 \]
\[
 r\dot{\theta}(\theta_{2}) = r^{n} \sin{((n-1)\theta_{2})} \geq 0.
\]
Again, we consider the component
\[ \dot{z}=z^{n}+\mu P_{1}(t)z,\]
with $0 \leq \mu \leq 1$. These equations do not have singular periodic solutions. But $\dot{z}=z^{n}$ has $n$ periodic solutions.
This completes the proof of the part about the number of complex periodic solutions.

Next, we prove the part about the number of real periodic solutions. The condition on $P_{1}$ implies that either
$\dot{r} > 0$ or $\dot{r}<0$ in the disk $r \leq 1$. Hence, periodic solutions do not lie in $r \leq 1$.
It follows from the phase portrait
that periodic solutions lie in $r >1 $. Hence a real periodic solution, $x(t)$,
lies in the intervals $(-\infty , -1)$ and $(1, \infty)$.

Consider the real differential equation (2.1).
We show that the second derivative of the displacement function $q''(c)$ does not change sign in $c>1$ and in $c<-1$.
This implies that $q$ has at most two zeros in any of these intervals.

When $x>1$,
\[ f_{2} = n(n-1) x^{n-2}+(n-1)(n-2)x^{n-3}P_{n-1}+\cdots+2P_{2}\]
We use the upper bound  for $|P_{i}|\leq k = \frac{n}{(n-2)^{2}}$. It is clear that this bound is larger than the upper bound
given in the statement of the theorem.
\[ f_{2} \geq (n-1)x^{n-3}(nx -(n-2)^{2}k) \geq \frac{(n-1)(n-4)}{n-3} > 0.\]
If $n$ is even and $x<-1$ then
\[f_{2} \geq -(n-1)x^{n-3}(-nx -(n-2)^{2}k) \geq (n-1)(n-(n-2)^{2}k)\geq \frac{(n-1)(n-4)}{n-3} > 0.\]
If $n$ is odd and $x<-1$ then
\[ f_{2} \leq (n-1)x^{n-3}(nx +(n-2)^{2}k) \leq \frac{(n-1)(4-n)}{n-3}  < 0.\]

In the case $n$ is odd, $q(c)$ has the sign of $c$ when $|c|$ is large.
Therefore, the number of real periodic solutions is odd.
\end{proof}

\begin{proof}
\emph{(Theorem 1.2:)}\\
As we explained in the proof of Theorem 1.1, we may assume that $K=1$.
We follow the same steps of proving Theorem 1.1 with some appropriate modifications. In the current case, if  $1<r<\rho$ then
\[
r\dot{\theta}(\theta_{1}) \leq (- r^{n}+P_{n-2}r^{n-2}+|P_{0}|)\sin{(\frac{k\pi}{n-2})}+(n-4)Cr^{n-3},
\]
where $C$ is the upper bound for $|P_{i}(t)|$ for $i \neq 0,1,n-2$.
Hence, $r\dot{\theta}(\theta_{1}) < 0$ when $r>1$ and $ (n-4)C \leq (1-|P_{0}|-P_{n-2})\sin{(\frac{\pi}{n-2})}$.
Similarly, we have
\[
r\dot{\theta}(\theta_{2}) \geq (r^{n}-P_{n-2}r^{n-2}-|P_{0}|)
\sin{(\frac{(k-1)\pi}{n-2})}-(n-4)Cr^{n-3}.
\]
When $r =1$, we have
\[
\dot{r}(r=1) \leq 1-P_{n-2} +(n-3)C+P_{1}+|P_{0}| \leq 0,
\]
when $P_{1} \leq -1-|P_{0}| - (n-3)C+P_{n-2}$. In this case, no solution leaves the disk $r<1$.

Similarly, if $P_{1} \geq 1+|P_{0}| + (n-3)C-P_{n-2}$, then $\dot{r}(r=1)\geq 0$.

Under the hypotheses in the statement of Theorem 1.2, no singular periodic solution can enter $D$ from $G_{k}$ and leave $D$,
when $k$ is odd and $k \leq n-1$. Since the phase portrait is symmetric about the $x-$axis, if $k$ is odd and $n \leq k \leq 2n-3$,
no singular periodic solution enters $D$ from $G_{k}$ and leaves $D$. It follows that the equation does not have a singular
periodic solution.

Now, consider the class of equations
\[
\dot{z}=z^{n}+\lambda[P_{n-3}(t)z^{n-3}+\cdots+P_{2}(t)z^{2}+P_{0}(t)]+P_{n-1}(t)z^{n-1}+P_{n-2}(t)z^{n-2}+P_{1}(t)z,
\]
with $0 \leq \lambda \leq 1$. For any equation in this family, the conditions in the statement of Theorem 1.2 are satisfied
and hence any equation does not have
a singular periodic solution. Since the equations in this family are in the same component of $\mathcal{A}$,
it follows from Lemma 2.1, that these equations have the same number of periodic solutions.
The equation
\[ \dot{z}=z^{n}+P_{n-1}(t)z^{n-1}+P_{n-2}(t)z^{n-2}+P_{1}(t)z \]
belongs to this family. This equation has $n$ periodic solutions. The result was proved in \cite{a1} for the case that
$P_{1}(t) \equiv 0$.
However the proof of this result involves only $\dot{\theta}$ and does not depend on $P_{1}(t)$. Therefore, this  result of
\cite{a1} is valid when
we add a linear term $P_{1}(t)z$. Now the result follows from Lemma 2.1.
\end{proof}
\begin{proof}
\emph{(Theorem 1.3:)}\\
We consider the real differential equation (2.1). We show that the $j$th derivative of $q(c)$ does not change sign in $c>1$
and in $c<-1$, and
the $(j-1)$th derivative of $q(c)$ does not change sign in $|c|<1$. If the $i$th derivative of a function does not change sign
on an interval then the
function has at most $i$ zeros in this interval. This implies that $q$ has at most $3j-1$ zeros. Now, we give the proof of
each case.
\begin{itemize}
\item[(i)]
If $x>1$ then
\[f_{1} \geq n x^{n-1} -(n-1)^{2} \frac{n}{(n-1)^{2}} x^{n-2} \geq nx^{n-2}(x-1) \geq 0.\]
If $x <-1$ and $n$ is even then
\[f_{1} \leq n x^{n-1} +(n-1)^{2} \frac{n}{(n-1)^{2}} x^{n-2} \leq nx^{n-2}(x+1) \leq 0.\]
If $x <-1$ and $n$ is odd then
\[f_{1} \geq n x^{n-1} +(n-1)^{2} \frac{n}{(n-1)^{2}} x^{n-2} \geq nx^{n-1}(1+\frac{1}{x}) \geq 0.\]
If $0<x<1$ then
\[f \leq x^{n-1}(x-\frac{n}{n-1}) + P_{0} \leq \frac{-1}{n-1}+P_{0}.\]
If $-1<x<0$ and $n$ is even then
\[f \leq -x^{n-1}(\frac{n}{n-1}-x) + P_{0} \leq \frac{2n-1}{n-1}+P_{0}.\]
If $-1<x<0$ and $n$ is odd then
\[f \leq x^{n-1}(x+\frac{n}{n-1}) + P_{0} \leq \frac{1}{n-1}+P_{0}.\]
Therefor, $f_{1} \geq 0$ in $x>1$ and $f_{1} \leq 0$ in $x<-1$. Moreover, $f \leq 0$ if $P_{0} \leq -\frac{2n-1}{n-1}$.
Similarly, $f \geq 0$ if $P_{0} \geq \frac{2n-1}{n-1}$.
\item[(ii)]
If $x>1$ then
\[f_{2} \geq n(n-1) x^{n-2} -(n-1)(n-2)^{2} \frac{n}{(n-2)^{2}} x^{n-2} \geq n(n-1)x^{n-2}(x-1) \geq 0.\]
If $x <-1$ and $n$ is even then
\[f_{2} \geq n(n-1) x^{n-2} +(n-1)^{2} \frac{n}{(n-1)^{2}} x^{n-3} \geq nx^{n-2}(1+\frac{1}{x}) \geq 0.\]
If $x <-1$ and $n$ is odd then
\[f_{2} \leq n(n-1) x^{n-2} +(n-1)^{2} \frac{n}{(n-1)^{2}} x^{n-3} \leq nx^{n-3}(x+1) \leq 0.\]
If $0<x<1$ then
\[f_{1}\leq nx^{n-2}(x+\frac{n-1}{n-2}) + P_{1} \leq \frac{n(2n-3)}{n-2}+P_{1}.\]
If $-1<x<0$ and $n$ is even then
\[f_{1} \leq nx^{n-2}(x+\frac{n-1}{n-2}) + P_{1} \leq \frac{n(2n-3)}{n-2}+P_{1}.\]
If $-1<x<0$ and $n$ is odd then
\[f_{1} \leq nx^{n-2}(x-\frac{n-1}{n-1}) + P_{1} \leq \frac{n(2n-3)}{n-2}+P_{1}.\]
Therefor, $f_{2} \geq 0$ in $x>1$ and $f_{1} \leq 0$ in $x<-1$. Moreover, $f \leq 0$ if $P_{1} \leq -\frac{n(2n-3)}{n-2}$.
Similarly, $f_{1} \geq 0$ if $P_{1} \geq \frac{n(2n-3)}{n-2}$.
\item[(iii)]
If $x>1$ then
\[f_{3} \leq n(n-1)(n-2) x^{n-3} -(n-4)^{2} \frac{n}{(n-3)^{2}} x^{n-4} \leq n(n-1)(n-2)x^{n-4}(x-1) \leq 0.\]
If $x <-1$ and $n$ is even then
\[f_{3} \leq n(n-1)(n-2) x^{n-3} -(n-4)^{2} \frac{n}{(n-3)^{2}} x^{n-4} \leq n(n-1)(n-2)x^{n-4}(x+1) \leq 0.\]
If $x <-1$ and $n$ is odd then
\[f_{3} \geq n(n-1)(n-2) x^{n-3} +(n-4)^{2} \frac{n}{(n-3)^{2}} x^{n-4} \geq n(n-1)(n-2)x^{n-3}(1+\frac{1}{x}) \geq 0.\]
If $0<x<1$ then
\[f_{2} \leq n(n-1)x^{n-3}(x+\frac{n-2}{n-3}) + P_{2} \leq \frac{n(n-1)(2n-5)}{n-3}+P_{2}.\]
If $-1<x<0$ and $n$ is even then
\[f_{2} \leq n(n-1)x^{n-3}(x-\frac{n-2}{n-3}) + P_{2} \leq \frac{n(n-1)(2n-5)}{n-3}+P_{2}.\]
If $-1<x<0$ and $n$ is odd then
\[f_{2} \leq n(n-1)x^{n-3}(x+\frac{n-2}{n-3}) + P_{2} \leq \frac{n(n-1)}{2n-5}+P_{2}.\]
Therefor, $f_{3} \geq 0$ in $x>1$ and $f_{3} \leq 0$ in $x<-1$. Moreover, $f_{2} \leq 0$ if $P_{2} \leq -\frac{n(n-1)(2n-5)}{n-3}$.
Similarly, $f_{2} \geq 0$ if $P_{2} \geq \frac{n(n-1)(2n-5)}{n-3}$.
\end{itemize}
\end{proof}
\begin{rema}
\begin{itemize}
\item[(i)]
Consider the differential equation
\begin{equation}
\label{e3.2}
\dot{z}=P_{n}(t) z^{n}+P_{n-1}(t)z^{n-1}+\dots+ P_{1}(t)z+P_{0}(t)
\end{equation}
If the leading coefficient, $P_{n}(t)$, does not vanish anywhere then the transformation of
the independent variable
\[
t \mapsto \int_{0}^{t}P_{n}(s) ds
\]
reduces the equation into a similar equation but with a leading coefficient equals 1.
\[ \dot{z}= z^{n}+\frac{P_{n-1}(t)}{P_{n}(t)} z^{n-1}+\dots+ \frac{P_{1}(t)}{P_{n}(t)} z+\frac{P_{0}(t)}{P_{n}(t)} .\]
Similar results can be proved for equation (3.2) with appropriate modifications.
\item[(ii)]
The conditions in our results on $|P_{i}|$ can be replaced by conditions on $\sum_{}^{}|P_{i}|$.
\end{itemize}
\end{rema}
\section{The Number of Limit Cycles}
We consider the two-dimensional autonomous system
\begin{equation}\label{e4.1}
\begin{cases}
\dot{x}=\lambda x - y + x(R_{1}(x,y)+R_{2}(x,y)+\cdots+R_{n-1}(x,y))\\
\dot{y}=x+\lambda y + y(R_{1}(x,y)+R_{2}(x,y)+\cdots+R_{n-1}(x,y)),
\end{cases}
\end{equation}
where $R_{i}$ is a homogeneous polynomials of degree $i$ for $i=1,2,\cdots,n-1$, and $\lambda$ is a real number.
The system in polar coordinates becomes
\[
\dot{r}=r^{n}R_{n-1}(\cos{\theta},\sin{\theta})+r^{n-1}R_{n-2}(\cos{\theta},\sin{\theta})+\cdots+
r^{2}R_{1}(\cos{\theta},\sin{\theta})+\lambda r,\]
\[\dot{\theta}=1.\]

Limit cycles of system (4.1) correspond to $2\pi-$periodic solutions of
\[
\frac{dr}{d \theta}=R_{n-1}r^{n}+R_{n-2}r^{n-1}+\cdots+R_{1}r^{2}+\lambda r.
\]
For this equation, $q'(0) < 0$ if $\lambda < 0$. If $\lambda =0$ and $\int_{0}^{2\pi} R_{1}(\theta) d\theta <0$ then $q''(0)<0$.
The following result follows from these remarks and the main results.
\begin{coro}
Suppose that
\[0< m \leq |R_{n-1}| \leq M.\]
\begin{itemize}
\item[(i)]
If
\[ |R_{i}| \leq k =\frac{m}{n-3} \sin{  (\frac{\pi}{n-2}) } ,\,\, i \neq n-1\]
 and
\[
| \lambda | \geq M+ \frac{m(n-2)}{n-3}\sin{  (\frac{\pi}{n-2}) }
\]
then system (4.1) has at most two limit cycles.
\item[(ii)]
Let $H=\max{|R_{n-3}|}$. If
\[ R_{n-3} \leq 0,\]
\[ |R_{i}| \leq  \frac{m}{n-4} \sin{(\frac{\pi}{n-2})},\,\, i\neq n-1,n-3\]
 and
\[
| \lambda | \geq M+H+ \frac{m(n-3)}{n-4}\sin{(\frac{\pi}{n-2})}
\]
then system (4.1) has at most $n-1$ limit cycles.
\item[(iii)]
If
\[|R_{i}| < \frac{nm}{(n-2)^{2}},\,\,i \neq n-1,\]
and
\[|\lambda| > \frac{nM(2n-3)}{n-2} \]
then system (4.1) has at most three limit cycles.
\item[(iv)]
If
\[|R_{i}| < \frac{nm}{(n-3)^{2}},\,\,i \neq 1,n-1,\]
and
\[|R_{1}| > \frac{nM(n-1)(2n-5)}{2(n-3)} \]
then system (4.1) has at most five limit cycles.
\end{itemize}
In parts (i), (ii), and (iii), there is at least one limit cycle when $\lambda <0$. In part (iv), there is a limit cycle
when $\lambda <0$, or when
$\lambda =0$ and $\int_{0}^{2\pi} R_{1}(\theta) d \theta <0$. In each case the origin is stable and the limit cycle is unstable.
\end{coro}
\begin{rema}
Formulae for the fourth and fifth derivatives of $q$ are given in \cite{a4}. Their use is more subtle than that of the first
three derivatives. In particular, $\frac{\partial^{4}f}{\partial x^{4}} > 0$ does not imply that $q^{iv}$ does not change sign.
\end{rema}


\begin{thebibliography}{00}
     \bibitem{a5} M. A. M. Alwash; {\em Complex centers of polynnomial differential equations},
     Electron. J. Differential Equations, 101 (2007), 1-15.
     \bibitem{a1} M. A. M. Alwash; {\em Periodic solutions of Abel differential equations}, J. Math. Anal. Appl.,
     329 (2007), 1161-1169
     \bibitem{a2} M. A. M. Alwash; {\em Periodic solutions of polynnomial non-autonous differential equations},
     Electron. J. Differential Equations, 84 (2005), 1-8.
    \bibitem{a3} M. A. M. Alwash; {\em Periodic solutions of a quartic
    differential equation and Groebner bases}, J. of Comp. and Appl. Math., 75 (1996), 67-76.
    \bibitem{a4} M. A. M. Alwash and N. G. Lloyd; {\em Periodic solutions of a quartic non-autonomus equation},
    Nonlinear Analusis, 11 (1987), 809-820.
    \bibitem{c1} M. Calanchi and B. Ruf; {\em On the number of closed solutions for polynomial ODE's
    and a special of Hilbert's 16th problem}, Adv. Differential Equations, 7 (2002), 197-216.
    \bibitem{c2} M. Calanchi and B. Ruf; {\em Hilbert type numbers for polynomial ODE's},
    Progress Nonlinear Differential Equations: Methods, Models, and Applications, 54(2003), 53-60.
    \bibitem{g1} A. Gasull and A. Guillamon; {\em Limit cycles for generalizes Abel equations}, Internat.
    J. Bifur. Chaos Appl. Sci. Engrg., 16 (2006), 3737-3745.
    \bibitem{i1} Y. Ilyashenko; {\em Hilbert-type numbers for Abel equations, growth and zeros of holomrphic functions},
    Nonlinearity, 13 (2000), 1337-1342.
    \bibitem{l1} A. Lins Neto; {\em On the number of solutions of the equations}
    $\frac{dx}{dt}=\sum_{j=0}^{n}a_{j}(t)x^{j}, 0 \leq t \leq 1,$ {\em for which} $x(0)=x(1)$, Inv. Math. 59 (1980), 67-76.
    \bibitem{l2} N.G. Lloyd; {\em The number of periodic solutions of the equation}
    $\dot{z}=z^{N}+p_{1}(t)z^{N-1}+\dots+p_{N}(t)$, Proc. London Math. Soc. (3) 27 (1973), 667-700.
    \bibitem{l3} N.G. Lloyd; {\em A note on the number of limit cycles in certain two-dimensional systems},
    J. London Math. Soc. 20 (1979), 277-286.
    \bibitem{p1} A. A. Panov; {\em The number of periodic solutions of polynomial differential
    equations}, Math. Notes 64 (1998), 622-628.
    \bibitem{p2} A. A. Panov; {\em On the diversity of Poincare mappings for cubic equations with variable coefficients},
    Functonal Analysis and its Applications 33 (1999), 310-312.
    \bibitem{p3} V. A. Pliss; {\em Nonlocal problems in the theory of oscillations}, Academic Press, New York (1966).
    \bibitem{s1} S. Smale; {\em Mathematical problems for the next century}, Mathematics: Frontiers and Perspectives,
    AMS (2000), 271-294.
\end{thebibliography}
\end{document}